\documentclass[fleqn]{mat01}
\usepackage{times,mathtimy,amssymb,latexsym}
\begin{document}

\setcounter{page}{197} \firstpage{197}

\title{Some further remarks on good sets}

\markboth{K Gowri Navada}{Some further remarks on good sets}

\author{K GOWRI NAVADA}

\address{Department of Mathematics, Periyar University, Salem~636~011, India\\
\noindent E-mail: gnavada@yahoo.com}

\volume{117}

\mon{May}

\parts{2}

\pubyear{2007}

\Date{MS received 24 January 2006; revised 19 August 2006}

\begin{abstract}
We show that in $n$-fold cartesian product, $n\geq 4$, a related
component need not be a full component. We also prove that when
$n\geq 4$, uniform boundedness of lengths of geodesics is not a
necessary condition for boundedness of solutions of $(1)$ for
bounded function $f$.
\end{abstract}


\maketitle

\section{Introduction and preliminaries}

The purpose of this note is to answer two questions about good
sets raised in \cite{3} and \cite{4}.

Let $X_{1},X_{2},\dots ,X_{n}$ be nonempty sets and let $\Omega
=X_{1}\times X_{2}\times \dots \times X_{n}$ be their cartesian
product. We will write $\vec{x}$ to denote a point
$(x_{1},x_{2},\dots ,x_{n})\in \Omega$. For each $1\leq i\leq n,
\Pi_{i}$ denotes the canonical projection of $\Omega$ onto
$X_{i}$.

A subset $S\subset \Omega$ is said to be {\it good}, if every
complex-valued function $f$ on $S$ is of the form:
\begin{equation}
f(x_{1},x_{2},\dots ,x_{n})=u_{1} (x_{1}) +u_{2} (x_{2})+\cdots
+u_{n} (x_{n} ),\ (x_{1}, x_{2},\dots, x_{n})\in S,
\end{equation}
for suitable functions $u_{1},u_{2},\dots , u_{n}$ on
$X_{1},X_{2},\dots ,X_{n}$ respectively (\cite{3}, p.~181).

For a good set $S$, a subset $B\subset \bigcup_{i=1}^{n} \Pi_{i}S$
is said to be a {\it boundary set of} $S$, if for any
complex-valued function $U$ on $B$ and for any $f\hbox{\rm :}\
S\longrightarrow \mathbb{C}$, eq.~(1) subject to
\begin{align*}
u_{i}\vert_{B\cap \Pi_{i}S}=U\vert_{B\cap \Pi_{i}S},\quad 1\leq
i\leq n,
\end{align*}
admits a unique solution. For a good set there always exists a
boundary set (\cite{3}, p.~187).

A subset $S\subset \Omega$ is said to be {\it full}, if $S$ is a
maximal good set in $\Pi_{1}S\times \Pi_{2}S\times\dots
\times\Pi_{n}S$.

A set $S \subset \Omega$ is full if and only if it has a boundary
consisting of $n-1$ points (\cite{3}, Theorem~3, page~185).

If a set $S$ is good, maximal full subsets of $S$ form a partition
of $S$. They are called {\it full components} of $S$ (\cite{3},
p.~183).

Two points $\vec{x},\vec{y}$ in a good set $S$ are said to be {\it
related}, denoted by $\vec{x} R \vec{y}$, if there exists a finite
subset of $S$ which is full and contains both $\vec{x}$ and
$\vec{y}$. $R$ is an equivalence relation, whose equivalence
classes are called {\it related components} of $S$. The related
components of $S$ are full subsets of $S$ (ref. \cite{3}).

\section{Example of a full set which is not related}

For a good set in two dimensions, full components are same as
related components and are called {\it linked components}
(\cite{2}, p.~60). In p.~190 of \cite{3}, the question whether
full components are the same as related components for $n\geq 3$,
is raised. Theorem~1 of \cite{5} answers the question partially,
where it is proved that a full set with finitely many related
components is itself related.

Here we prove that when the dimension $n$ is $\geq 4$, a full
component need not be a related component, by giving an example of
a full set with infinitely many related\break components.

Our example, in a four-fold cartesian product, will consist of
countable number of points $S=\{\vec{y_{1}},\vec{y_{2}},
\vec{y_{3}},\dots\}$ such that for each $n$ the subset $S_{n}$ of
first $n$ points of this set will be a good set and have a
boundary consisting of four or five points depending on whether
$n$ is even or odd. Moreover, any boundary of $S_{n}$ will
necessarily contain two or one point from the coordinates of
$\vec{y}_{n}$, depending on whether $n$ is odd or even, so that
eventually these points of the boundary disappear and $S$ will
have boundary with only three points. So $S$ will be full. On the
other hand, $S$ will have related components consisting of
single\break points.

Let
$\{x_{1},x_{2},x_{3},x_{4},\alpha_{1},\alpha_{2},\alpha_{3},\dots,\alpha_{n},\dots\}$
be a set of distinct symbols. The $j$th term of this sequence will
be called the $j$th symbol. Then $x_{i}$ will be $i$th symbol,
$1\leq i\leq 4$, while $\alpha_{j}$ will be the $(j+4)$th symbol
for $j\geq 1$.

The countable infinite set $S$ in four dimensions is defined as
$S=\{\vec{y_{1}}, \vec{y_{2}}, \vec{y_{3}},\dots \}$,\break where
\begin{align*}
\vec{y_{1}} &= (x_{1},x_{2},x_{3},x_{4}),\\[.3pc]
\vec{y_{2}} &=
(x_{1},x_{2},\alpha_{1},\alpha_{2}),\\[.3pc]
\vec{y_{3}} &=(\alpha_{3},\alpha_{4},x_{3},\alpha
_{2}),\\[.3pc]
\vec{y_{4}}
&=(\alpha_{3},\alpha_{4},\alpha_{1},x_{4}),\\[.3pc]
\vec{y_{5}}
&=(x_{1},\alpha_{4},\alpha_{5},\alpha_{6}),\\[.3pc]
\vec{y_{6}} &=(\alpha_{3},x_{2},\alpha_{5},\alpha_{6} ),
\end{align*}
and so on. In general, for $n\geq 2$,
\begin{align*}
\vec{y}_{4n-3} &= (x_{1} ,\alpha_{4n-4},\alpha_{4n-3}
,\alpha_{4n-2}),\quad \vec{y}_{4n-2}
=(\alpha_{4n-5},x_{2},\alpha_{4n-3},\alpha_{4n-2}),\\[.3pc]
\vec{y}_{4n-1} &= (\alpha_{4n-1},\alpha_{4n},x_{3} ,\alpha_{4n-2}
),\quad \vec{y}_{4n} =(\alpha_{4n-1},\alpha_{4n} ,\alpha_{4n-3}
,x_{4} ).
\end{align*}

We prove that this set is full and singletons are its related
components.

Let $S_{n}=\{\vec{y_{1}},\vec{y_{2}},\dots,\vec{y_{n}}\}.$
Consider the matrix $M_{n}$ corresponding to the set $S_{n}$,
called the {\it matrix} of $S_{n}$ (p.~58 of \cite{2}), whose rows
correspond to the points $\vec{y_{i}}$, $1\leq i\leq n$, and whose
columns correspond to the symbols occurring in the points of
$S_{n}$. The $ij$th entry of this matrix is $1$, if the $j$th
symbol occurs in the point $\vec{y_{i}}$. Otherwise, the $ij$th
entry\break is $0$.

If $n$ is odd, say $n=2m-1$, then $M_{2m-1}$ consists of $2m-1$
rows and $2m+4$ columns corresponding to the symbols
$\{x_{1},x_{2},x_{3},x_{4},\alpha_{1},\alpha_{2},\dots
,\alpha_{2m-1},\alpha_{2m}\}.$ If $n$ is even, $n=2m $, then there
are $2m$ rows and $2m+4$ columns corresponding to
$\{x_{1},x_{2},x_{3},x_{4},$ $\alpha _{1},\alpha_{2},\dots
,\alpha_{2m-1},\alpha _{2m}\}$ in $M_{2m}$. The matrix of $S$ is
defined similarly.\pagebreak

We will prove by induction that for each $n,$ $S_{2n}$ is a good
set and has boundary $\{x_{1},x_{2},x_{3},$ $\alpha_{2n}\}$ (or
$\{x_{1},x_{2},x_{3},\alpha _{2n-1}\}$). Clearly the statement
holds for $n=1.$ Assume that the statement holds for $n=m-1$. We
have to show that the statement holds for $n=m.$

Since $S_{2m-2}$ is good and has a boundary consisting of
$\{x_{1},x_{2},x_{3},\alpha _{2m-2}\}$ it is clear that $S_{2m-1}$
is also good and has a boundary consisting of
$\{x_{1},x_{2},x_{3},\alpha_{2m-1},\alpha_{2m}\}.$ (Note that
$\vec{y}_{2m-1}$ has two new coordinates $\alpha _{2m-1}$ and
$\alpha _{2m}$ not occurring in $\vec{y}_{1},\dots
,\vec{y}_{2m-2}$). For any $f$ on $S_{2m-1},$ a solution
$u_{1},u_{2},u_{3},u_{4}$ satisfying
\begin{align*}
f(z_{1},z_{2},z_{3},z_{4}) &=
u_{1}(z_{1})+u_{2}(z_{2})+u_{3}(z_{3})+u_{4}(z_{4}),\\[.3pc]
&\quad\, (z_{1},z_{2},z_{3},z_{4})\in S_{2m-1},
\end{align*}
is uniquely determined once we fix the values of $u_{i},\ 1\leq
i\leq 4$, on the boundary points.

We drop the columns corresponding to the symbols
$x_{1},x_{2},x_{3},\alpha_{2m-1},\alpha_{2m}$ from $M_{2m-1}$ and
get an invertible $(2m-1)\times (2m-1)$-matrix $N_{2m-1}$ with
entries zeros and ones given below:
\begin{equation*}
N_{2m-1}=
\left[\begin{array}{c@{\quad}c@{\quad}c@{\quad}c@{\quad}c@{\quad}c@{\quad}c@{\quad}c@{\quad}c@{\quad}c}
1 & 0 & 0 & 0 & 0 & 0 & . & . & . & 0 \\[.2pc]
0 & 1 & 1 & 0 & 0 & 0 & . & . & . & 0 \\[.2pc]
0 & 0 & 1 & 1 & 1 & 0 & 0 & . & . & 0 \\[.2pc]
1 & 1 & 0 & 1 & 1 & 0 & 0 & . & . & 0 \\[.2pc]
0 & 0 & 0 & 0 & 1 & 1 & 1 & . & . & 0 \\[.2pc]
0 & 0 & 0 & 1 & 0 & 1 & 1 & . & . & 0 \\
. & . & . & . & . & . & . & . & . & . \\[.2pc]
. & 0 & 0 & 0 & 0 & 0 & 0 & 1 & 1 & 1 \\[.2pc]
. & 0 & 0 & 0 & 0 & 0 & 1 & 0 & 1 & 1 \\[.2pc]
. & 0 & 0 & 0 & 0 & 0 & 0 & 0 & 0 & 1
\end{array} \right].
\end{equation*}
To show that $S_{2m}$ is good and has a boundary consisting of
$\{x_{1},x_{2},x_{3},\alpha_{2m}\}$, consider the matrix of
$S_{2m}$.

If we drop from the matrix $M_{2m}$ the columns corresponding to
$x_{1},x_{2},x_{3},\alpha_{2m},$ we get a $2m\times 2m$-matrix
given below.
\begin{equation*}
N_{2m}=
\left[\begin{array}{c@{\quad}c@{\quad}c@{\quad}c@{\quad}c@{\quad}c@{\quad}c@{\quad}c@{\quad}c@{\quad}c}
1 & 0 & 0 & 0 & 0 & 0 & . & . & . & 0 \\[.2pc]
0 & 1 & 1 & 0 & 0 & 0 & . & . & . & 0 \\[.2pc]
0 & 0 & 1 & 1 & 1 & 0 & 0 & . & . & 0 \\[.2pc]
1 & 1 & 0 & 1 & 1 & 0 & 0 & . & . & 0 \\[.2pc]
0 & 0 & 0 & 0 & 1 & 1 & 1 & . & . & 0 \\[.2pc]
0 & 0 & 0 & 1 & 0 & 1 & 1 & . & . & 0 \\
. & . & . & . & . & . & . & . & . & . \\[.2pc]
. & 0 & . & . & . & 0 & 1 & 1 & 1 & 0 \\[.2pc]
. & 0 & . & . & . & 1 & 0 & 1 & 1 & 0 \\[.2pc]
. & 0 & . & . & . & 0 & 0 & 0 & 1 & 1 \\[.2pc]
. & 0 & 0 & 0 & 0 & 0 & 0 & 1 & 0 & 1
\end{array}\right].
\end{equation*}
Its initial $(2m-1)\times (2m-1)$ matrix is the matrix $N_{2m-1}$
obtained above which is invertible, and hence has its rows
linearly independent. It is clear from this that rows of $N_{2m}$
are also linearly independent, so that $N_{2m}$ is invertible.
This proves that $S_{2m}$ is good and
$\{x_{1},x_{2},x_{3},\alpha_{2m}\}$ is its boundary. We can see
similarly that $\{x_{1},x_{2},x_{3},\alpha_{2m-1}\}$ is also a
boundary for $S_{2m}$.

We now prove that $\{x_{1},x_{2},x_{3},\alpha_{2m}\}$ (or
$\{x_{1},x_{2},x_{3},\alpha_{2m-1}\}$) is not a boundary of
$S_{2n}$ for any $m<n$. Indeed, if $\{x_{1},x_{2},x_{3},
\alpha_{2m}\}$ is a boundary of $S_{2n}$ for some $m<n,$ then for
any $f$ on $S_{2m+2}$ there is a solution
$u_{1},u_{2},u_{3},u_{4}$ of
\begin{align*}
f(z_{1},z_{2},z_{3},z_{4}) &= u_{1}(z_{1})+u_{2}
(z_{2})+u_{3}(z_{3})+u_{4}(z_{4}),\\[.3pc]
&\quad\, (z_{1},z_{2},z_{3},z_{4} )\in S_{2m+2},
\end{align*}
where values of $u_{1},u_{2},u_{3}$ and an appropriate $u_i$ are
preassigned on $x_{1},x_{2},x_{3},\alpha_{2m}$ respectively.

Since $\{x_{1},x_{2},x_{3},\alpha_{2m}\}$ is a boundary of
$S_{2m}$, $u_{4}(x_{4})$ and the value of appropriate $u_{i}$ is
determined on $\alpha_{2m-1}$ by the values of the function on the
points of $S_{2m}$.

Now if $2m\equiv 0({\rm mod}\ 4)$, then
\begin{align*}
f( \vec{y}_{2m+1})&=u_{1} (x_{1} )+u_{2}
(\alpha_{2m})+u_{3} (\alpha_{2m+1})+u_{4} (\alpha_{2m+2}),\\[.3pc]
f(\vec{y}_{2m+2})&=u_{1}(\alpha_{2m-1})+u_{2}(x_{2})+u_{3}(\alpha_{2m+1})+u_{4}(\alpha_{2m+2}),
\end{align*}
which clearly do not hold together if
\begin{equation*}
f(\vec{y}_{2m+2})\neq
f(\vec{y}_{2m+1})-u_{1}(x_{1})-u_{2}(\alpha_{2m})+u_{1}
(\alpha_{2m-1})+u_{2}(x_{2}).
\end{equation*}
The case $2m\equiv 2({\rm mod}\ 4)$ can be treated similarly and
we see that the set $\{x_{1},x_{2},x_{3},\alpha_{2m}\}$ cannot be
a boundary of $S_{2n}$ for $n>m$. The case of
$\{x_{1},x_{2},x_{3},\alpha_{2m-1}\}$ is si\-mi\-lar. The set $S$
is good because every finite subset of it is good. To show that it
is full it is enough to show that any boundary $B$ of $S$ consists
of three points. Take a function $f$ on $S$ and fix the values of
$u_{i}$ on $\Pi_{i}S\cap B,\ 1\leq i\leq 4$. Then there exists a
unique solution $u_{i}$ on $\Pi_{i}S,\ i=1,2,3,4$ of $(1)$. This
gives a solution of $(1)$ on $S_{2n}$ with $f=f|_{S_{2n}}$, for
any $n\in \mathbb{Z}_{+}$. Fix values of $u_{i}$ on $B\cap
\Pi_{i}S_{2n},1\leq i\leq 4$. Then there is a solution of $(1)$
with $f=f|_{S_{2n}},$ and $u_{i}$ prescribed on $B\cap
\Pi_{i}S_{2n},\ 1\leq i\leq 4$, as above. If $|B|
>3 $, then for $n$ large enough we get $|B\cap (\cup
_{i=1}^{4}\Pi_{i}S_{2n}) | \geq 4$, and this will give a boundary
for $S_{2n}$ with, either more than four points or, a four-point
boundary which does not contain $\alpha_{2n-1}$ and $\alpha_{2n}$.
But this is not possible. This shows that $|B| =3$ which means
$S$\break is full.

No finite subset $K$ of $S,$ other than singleton, is full. We
prove this by showing that there is a point
$\vec{y}=(y_{1},y_{2},y_{3},y_{4})\notin K$ with $y_{i}\in \Pi
_{i}K$ for $ i=1,2,3,4$ such that $K\cup \{\vec{y}\}$ is also
good. Any two-point subset of $S$ is not full as any two points of
$S$ differ in at least two coordinates. So we can assume that
$|K|\geq 3$. Let $n$ be the least integer such that $K \subseteq
S_{2n}$. Then either $\vec{y}_{2n-1}$ or $\vec{y}_{2n}\in K$. If
$\vec{y}_{2n} \in K$, then let $\vec{z}= \vec{y}_{2n}$; otherwise
let $\vec{z}=\vec{y}_{2n-1}.$ Then $\alpha_{2n-1}$ and
$\alpha_{2n}$ are the coordinates of $\vec{z}$.

Let $\vec{y}$ be same as $\vec{z}$ except that the coordinate $
\alpha _{2n}$ is replaced by some other corresponding coordinate
of a point in $K$. Then $\vec{y}\notin K.$ Further $K\cup
\{\vec{y}\}$ is good: for which it is enough to show that
$S_{2n}\cup \{\vec{y}\}$ is good.

Consider the $(2n+1)\times (2n+1)$-matrix whose columns correspond
to the symbols $\{x_{4},\alpha_{1},\dots
,\alpha_{2n-1},\alpha_{2n}\}$ and rows correspond to the points of
$\{\vec{y}_{1},\dots ,\vec{y}_{2n},\vec{y}\}$.\pagebreak
\begin{equation*}
\left[\begin{array}{c@{\quad}c@{\quad}c@{\quad}c@{\quad}c@{\quad}c@{\quad}c@{\quad}c@{\quad}c@{\quad}c@{\quad}c}
1 & 0 & 0 & 0 & 0 & 0 & . & . & . & . & 0 \\[-.05pc]
0 & 1 & 1 & 0 & 0 & 0 & . & . & . & . & 0 \\[-.05pc]
0 & 0 & 1 & 1 & 1 & 0 & 0 & . & . & . & 0 \\[-.05pc]
1 & 1 & 0 & 1 & 1 & 0 & 0 & . & . & . & 0 \\[-.05pc]
0 & 0 & 0 & 0 & 1 & 1 & 1 & 0 & . & . & 0 \\[-.05pc]
0 & 0 & 0 & 1 & 0 & 1 & 1 & 0 & 0 & . & 0 \\[-.15pc]
. & . & . & . & . & . & . & . & . & . & . \\
. & 0 & 0 & 0 & 0 & 1 & 0 & 1 & 1 & 0 & 0 \\[-.05pc]
. & 0 & 0 & 0 & 0 & 0 & 0 & 0 & 1 & 1 & 1 \\[-.05pc]
. & 0 & 0 & 0 & 0 & 0 & 0 & 1 & 0 & 1 & 1 \\[-.05pc]
. & . & . & . & . & . & . & . & . & 1 & 0
\end{array}\right],
\end{equation*}
where the last row has a $1$ at the column corresponding to
$\alpha_{2n-1}$ and a $0$ in the column corresponding to
$\alpha_{2n} $. The first $2n$ rows of this matrix are linearly
independent as the set $S_{2n}=\{\vec{y}_{1},\dots ,
\vec{y}_{2n}\}$ is good. Further, it is clear that the last row is
not a linear combination of the first $2n$ rows. Therefore
$S_{2n}\cup \{\vec{y}\}$ is good. This proves that the related
components of $S$ are singletons.\vspace{-.2pc}

\section{A set with unbounded geodesic length which has bounded solutions for bounded functions $\pmb{f}$}

A modification of this example can be used to answer another
question.\vspace{-.05pc}

Let $S$ be a related set in $n$-fold cartesian product
$X_{1}\times X_{2}\times \dots \times X_{n}$.\vspace{-.05pc}

Given two points\vspace{-.05pc} $\vec{x},\vec{y}\in S$ the
smallest related subset of $S$ containing $\vec{x}$ and $\vec{y}$
is called the {\it geodesic} between $\vec{x}$\vspace{-.05pc} and
$\vec{y}$ (\cite{3}, p.~190). In \cite{4}, it is shown that if
there is an upper bound for lengths of geodesics in $S$, then the
solutions $u_{1},\dots ,u_{n}$ of $(1)$ are bounded whenever $f$
is bounded. Further, the author surmises that the converse may be
true. Namely, that boundedness of geodesic lengths is a necessary
condition for bounded solutions $u_{1},\dots ,u_{n},$ for bounded
$f$. This is true for $n=2$. However, for $n\geq 4$ we prove that
this is not a necessary condition.

Consider the set $S_{4n}\cup \{\vec{z}\}$\vspace{-.05pc} where the
point $\vec{z}=(\alpha_{4n-1}, x_{2},\alpha_{4n-3},x_{4}).$ We
prove that $S_{4n}\cup \{\vec{z}\}$ is full.\vspace{-.05pc} First
we note that the infinite matrix with rows corresponding to
$\{\vec{y}_{2},\vec{y}_{3},\dots \}$ and columns corresponding to
$\{\alpha_{1},\alpha_{2},\dots \}$ has an inverse given by
\begin{equation*}
\left[\begin{array}{cccccccccccc} \frac{1}{2} & -\frac{1}{2} &
\medskip \frac{1}{2} & 0 & 0 & 0 & .
& . & . & .& . & . \\[-.15pc]
\frac{1}{2} & \frac{1}{2} & -\frac{1}{2}\medskip & 0 & 0 &
0 & . & . & . & . & . & . \\[-.15pc]
-\frac{1}{4} & \frac{1}{4} & \medskip \frac{1}{4} & -\frac{1}{2} &
\frac{1}{2} & 0 & 0 & 0 & . & . & . & . \\[-.15pc]
-\frac{1}{4} & \frac{1}{4} & \medskip \frac{1}{4} & \frac{1}{2} &
-\frac{1}{2} & 0 & 0 & 0 & . & . & . & . \\[-.15pc]
\frac{1}{8} & -\frac{1}{8} & -\frac{1}{8}\medskip & \frac{1}{4} &
\frac{1}{4}& -\frac{1}{2} & \frac{1}{2} & 0 & 0 & 0 & . & . \\[-.15pc]
\frac{1}{8} & -\frac{1}{8} & \medskip -\frac{1}{8} & \frac{1}{4} &
\frac{1}{4} & \frac{1}{2} & -\frac{1}{2} & 0 & 0 & 0 & . & . \\[-.15pc]
-\frac{1}{16} & \frac{1}{16} & \medskip \frac{1}{16} &
-\frac{1}{8} & -\frac{1}{8} & \frac{1}{4} &
\frac{1}{4} & -\frac{1}{2} & \frac{1}{2} & 0 & 0 & . \\[-.15pc]
-\frac{1}{16} & \frac{1}{16} & \medskip \frac{1}{16} &
-\frac{1}{8} & -\frac{1}{8} & \frac{1}{4} & \frac{1}{4} &
\frac{1}{2} & -\frac{1}{2} & 0 & 0 & . \\[-.35pc]
. & . & \medskip . & . & . & . & . & . & . & . & . & . \\[-.15pc]
\frac{1}{2^{2n-1}} & -\frac{1}{2^{2n-1}} & \medskip
-\frac{1}{2^{2n-1}} & \frac{1}{2^{2n-2}} & . & . & . & . & . &
-\frac{1}{2} & \frac{1}{2} & . \\[-.15pc]
\frac{1}{2^{2n-1}} & -\frac{1}{2^{2n-1}} & \medskip
-\frac{1}{2^{2n-1}} & \frac{1}{2^{2n-2}} & . & . & . & . & . &
\frac{1}{2} & -\frac{1}{2} & . \\[-.35pc]
. & . & . & . & . & . & . & . & . & . & . & .
\end{array}\right].
\end{equation*}
Now consider the $4n\times 4n$-matrix $A_{4n}$ with rows
corresponding to the points $\{\vec{y_{2}}, \vec{y_{3}},\dots ,$
$\vec{y_{4n}},\vec{z}\}$\ and columns corresponding to
$\{\alpha_{1},\dots ,\alpha_{4n}\}$:
\begin{equation*}
A_{4n}= \left[
\begin{array}{c@{\quad}c@{\quad}c@{\quad}c@{\quad}c@{\quad}c@{\quad}c@{\quad}c@{\quad}c@{\quad}c@{\quad}c@{\quad}c}
1 & 1 & 0 & 0 & 0 & 0 & 0 & 0 & . & . & . & 0 \\[.2pc]
0 & 1 & 1 & 1 & 0 & 0 & 0 & 0 & . & . & . & 0 \\[.2pc]
1 & 0 & 1 & 1 & 0 & 0 & 0 & 0 & . & . & . & 0 \\[.2pc]
0 & 0 & 0 & 1 & 1 & 1 & 0 & 0 & . & . & . & 0 \\[.2pc]
0 & 0 & 1 & 0 & 1 & 1 & 0 & 0 & . & . & . & 0 \\[.2pc]
0 & 0 & 0 & 0 & 0 & 1 & 1 & 1 & . & . & . & 0 \\[.2pc]
0 & 0 & 0 & 0 & 1 & 0 & 1 & 1 & . & . & . & 0 \\
. & . & . & . & . & . & . & . & . & . & . & . \\[.2pc]
0 & 0 & 0 & . & . & . & . & . & 1 & 1 & 0 & 0 \\[.2pc]
0 & 0 & 0 & . & . & . & . & . & 0 & 1 & 1 & 1 \\[.2pc]
0 & 0 & 0 & 0 & 0 & . & . & . & 1 & 0 & 1 & 1 \\[.2pc]
0 & 0 & 0 & 0 & 0 & . & . & . & 1 & 0 & 1 & 0
\end{array}\right].
\end{equation*}

This matrix is invertible and the inverse is
\begin{equation*}
\begin{bmatrix}
\frac{1}{2}\medskip & -\frac{1}{2} & \frac{1}{2} & 0 & 0 & 0 & 0 & . & . & 0\\
\frac{1}{2}\medskip & \frac{1}{2} & -\frac{1}{2} & 0 & 0 & 0 & 0 & . & . & 0 \\
-\frac{1}{4}\medskip & \frac{1}{4} & \frac{1}{4} & -\frac{1}{2} &
\frac{1}{2} & 0 & 0 & . & . & 0 \\
-\frac{1}{4}\medskip & \frac{1}{4} & \frac{1}{4} & \frac{1}{2} &
-\frac{1}{2} & 0 & 0 & . & . & 0 \\
\frac{1}{8}\medskip & -\frac{1}{8} & -\frac{1}{8} & \frac{1}{4} &
\frac{1}{4} & -\frac{1}{2} & \frac{1}{2} & . & . & 0 \\
. & . & . & . & . & . & . & . & . & . \\
\frac{1}{2^{2n-1}}\medskip & -\frac{1}{2^{2n-1}} &
-\frac{1}{2^{2n-1}} & \frac{1}{2^{2n-2}} & . & . & . &
-\frac{1}{2} & \frac{1}{2} & 0 \\
\frac{1}{2^{2n-1}}\medskip & -\frac{1}{2^{2n-1}} &
-\frac{1}{2^{2n-1}} & \frac{1}{2^{2n-2}} & . & . & . & \frac{1}{2}
&-\frac{1}{2} & 0 \\
-\frac{1}{2^{2n-1}}\medskip & \frac{1}{2^{2n-1}} &
\frac{1}{2^{2n-1}} & -\frac{1}{2^{2n-2}} & . & . & . & \frac{1}{2}
&-\frac{1}{2} & 1 \\
0 & 0 & 0 & 0 & 0 & . & . & 0 & 1 & -1
\end{bmatrix}.
\end{equation*}
This proves that $S_{4n}\cup \{\vec{z}\}$ is good and
$\{x_{1},x_{2},x_{3}\}$ is its boundary. So this set is also full.
Note that the sums of the absolute values of entries in a row is
less than or equal to $3$, for any row, independent of $n$.

We prove that the geodesic between $\vec{y}_{1}$ and
$\vec{y}_{4n}$ in $ S_{4n}\cup \{\vec{z}\}$ is the whole set
$S_{4n}\cup \{\vec{z}\}.$ If $K\subset S_{4n}\cup \{\vec{z}\}$ is
a full set containing $\{\vec{y}_{1},\vec{y}_{4n}\}$ then we have
to show that $K=S_{4n}\cup \{\vec{z}\} $.

First we note that $\vec{z}\in K$. If $K\neq S_{4n}\cup
\{\vec{z}\}$, let $i$ be the number of points of $S_{4n}$ which
are not in $K$. These $i$ points should have at least $i$ symbols
occurring in them which are not in $\cup _{i=1}^{4}\Pi_{i}K$
because $K$ is full and, when we add these $i$ points to $K$ the
set remains good. These symbols are from $\{\alpha_{1},\dots
,\alpha_{4n-2}\}$ because $x_{1},x_{2},x_{3},x_{4},\alpha_{4n-1}$
and $\alpha _{4n}$ occur as co-ordinates in the points of $K$. Let
$\alpha_{j_{1}},\alpha_{j_{2}},\dots ,\alpha_{j_{i}}$ where
$j_{1}<j_{2}<\cdots <j_{i}$ be these symbols. If we prove that
these symbols are used by at least $i+1$ points of $S_{4n}$, we
get a contradiction because these $i+1$ points cannot be in $K$.

For this, we show that the $i$ columns of the matrix $A_{4n}$
corresponding to $\alpha_{j_{1}}, \alpha_{j_{2}},\dots
,\alpha_{j_{i}}$ have nonzero entries in at least $i+1$ rows of
$A_{4n} $. Let us first take the case when $\alpha_{j_{s}}\neq
\alpha_{1}$ or $\alpha _{2}$ for $s=1,2,\dots ,i$. Then these $i$
columns contain exactly $3i$ nonzero entries in them. Any row of
$A_{4n}$ contains at most three nonzero entries in it. But observe
that the row containing the last nonzero entry of the column
corresponding to $\alpha _{j_{i}}$ has only one nonzero entry in
these $i$ columns. So we need at least $i+1$ rows to cover all the
nonzero entries of these $i$ columns. Now assume that
$\alpha_{j_{s}}=\alpha_{1}$ (or $\alpha_{2})$ for some $s$. Then
the total number of nonzero entries in these $i$ columns is
$3(i-1)+2$. As in the previous case, there is a row containing
only one nonzero entry of these $i$ columns. So again it is easy
to see that we need at least $i+1$ rows to cover all the nonzero
entries in these $i$ columns. In the last case, when
$\alpha_{j_{s}}=\alpha_{1}$ and $\alpha_{j_{t}}=\alpha_{2}$ for
some $1\leq s,\ t\leq i$, the total number of nonzero entries in
these columns is $3(i-2)+2+2$. But in this case the first row
contains only two nonzero entries of these $i$ columns. As before
the row containing the last nonzero entry of the column
corresponding to $\alpha_{j_{i}}$ has only one nonzero entry of
these columns. So, again we need at least $i+1$ rows to cover the
nonzero entries of these columns. This contradiction proves that
$K=S_{4n}\cup \{\vec{z}\}$.

Since the bound on the absolute row sums of $A_{4n}^{-1}$ is
independent of $n$, it is clear how to construct a full set $S$ in
which geodesic lengths are not bounded, but solutions
$u_{1},\dots,u_{n}$ of $(1)$ are bounded for bounded $f$.

A~similar kind of construction is possible for the case $n=3$ also
and it will be communicated shortly.

\section*{Acknowledgement}

The author would like to thank Prof.~M~G~Nadkarni for fruitful
discussions and encouragement and for suggesting the problem. The
author also thanks IMSc for short visiting appointments which made
this work possible.

\end{document}